\documentclass{article}
\usepackage{amsmath}
\usepackage{amsfonts}
\usepackage{amssymb}
\usepackage{graphicx}

\input{tcilatex}
\begin{document}

\title{The Topological Structure of Question Theory}
\author{Shahid Nawaz\thanks{%
Email: sn439165@albany.edu} \\
{\small Department of Physics, University at Albany-SUNY, }\\
{\small Albany, NY 12222, USA.}}
\date{}
\maketitle

\begin{abstract}
A question is identified with a topology on a given set of irreducible
assertions. It is shown that there are three types of a question. Type-I
question generates sub-question, type-II question has a definite answer and
type-III question is irrelevant. We suggest that the most intelligent
machine asks type-II questions. We also claim that a truly intelligent
machine cannot be desireless. This work may prove useful in machine learning
and may open up new ways to understand mind.
\end{abstract}

\section{Introduction}

It has been pointed out by many authors that question is of fundamental
importance \cite{Cox, Cox79,Fry,Fry01, Fry02, Knuth, Knuth04, Knuth07,
Caticha}. In most cases we know the solution of a problem but we do not know
the right question that answers it. It also happens sometime that a problem
can be solved in many different ways. A truly intelligent machine is needed
to ask relevant questions. It is however unimportant, at this stage, that
some solutions are better than other.

A question is defined as the set of assertions that answers it. An older
definition of the question is that a question is a request for information.
The former definition was given by Richard T. Cox \cite{Cox, Cox79}, and
studied further by others \cite{Fry,Fry01, Fry02, Knuth, Knuth04, Knuth07,
Caticha}. We follow R.T. Cox and define a question is a topology on a given
set of irreducible assertions.

This paper is a self-contained development of the topological structure of
question theory. Some definitions and theorems are given in section-2 which
shall be used in the following sections. The paper is concluded in section-4.

\section{Some Definitions and Theorems}

The following definitions and theorems will be used in the sequel, see for
instant \cite{Munkres}.\newline
\textbf{Definition 2.1:\qquad }\textit{Let }$X$\textit{\ be a non-empty set.
A collection }$\mathcal{T}$ \textit{of subsets of }$X$\textit{\ is called a
topology on }$X$\textit{, if it meets the following requirements:}\newline
\textit{C1. }$\phi $\textit{\ and }$X$\textit{\ belong to }$\mathcal{T}$.%
\newline
\textit{C2. The union of any number of sets in }$\mathcal{T}$ \textit{belong
to} $\mathcal{T}$.\newline
\textit{C3. The intersection of any two sets (and hence of any finite number
of set}) in$\mathcal{\ T}$ \textit{belongs to} $\mathcal{T}$.\newline
\textbf{Definition 2.2:\qquad }\textit{Let }$\left( X,\mathcal{T}\right) $ 
\textit{be a topological space. A subset }$U$ \textit{of }$X$\textit{\ is
said to be open iff it belongs to }$\mathcal{T}$.\newline
The complement $X-U$ is called closed set and a set which is both open and
closed is called clopen.\newline
\textbf{Definition 2.3:\qquad }\textit{Let }$\left( X,\mathcal{T}\right) $%
\textit{\ be a topological space and let }$x$\textit{\ \ be a point of }$X$%
\textit{. \ A subset }$N$\textit{\ of }$X$\textit{\ is called a neighborhood
(nbhd in short) of} $x$ \textit{iff there exists an open set }$u$\textit{\
such that }$x\in u\subseteq N\subseteq X$\newline
\textbf{Definition 2.4:\qquad }\textit{Let }$x_{0}$ \textit{be a point in a
topological space }$\left( X,\mathcal{T}\right) $, \textit{then the set of
all nbhds} \textit{is called nhbd system of the point }$x_{0}$ \textit{and
is denoted by} $\mathcal{N}_{(x_{0})}$.\newline
\textbf{Theorem 2.1:\qquad }\textit{Let }$\left( X,\mathcal{T}\right) $ 
\textit{be a topological space and let }$\mathcal{N}_{\left( x_{0}\right) }$ 
\textit{be the nbhd system of the point }$x_{0}\in X$\textit{,} \textit{then}%
\newline
\textit{1.\qquad }$\mathcal{N}_{\left( x_{0}\right) }$ \textit{is non-empty.}%
\newline
\textit{2.\qquad The intersection of any two members of }$\mathcal{N}%
_{\left( x_{0}\right) }$ \textit{belongs to }$\mathcal{N}_{\left(
x_{0}\right) }$.\newline
\textit{3.\qquad If }$U$ \textit{is in }$\mathcal{N}_{\left( x_{0}\right) }$ 
\textit{and }$W$ \textit{is any set of X such that }$U\subset W$, \textit{%
then W is in }$\mathcal{N}_{\left( x_{0}\right) }$.\newline
\textbf{Definition 2.5:\qquad }\textit{Let }$X$ \textit{be a topological
space with topology }$\mathcal{T}$. \ \textit{If }$A$ \textit{is subset of }$%
X$\textit{,} \textit{the collection}%
\begin{equation*}
\mathcal{T}_{A}=\left\{ A\cap U:U\in \mathcal{T}\right\} 
\end{equation*}%
\textit{is a topology on }$A$, \textit{called the subspace topology,and }$A$ 
\textit{is called a subspace of }$X$.\newline
\textbf{Definition 2.6:\qquad }\textit{Let }$A$\textit{\ and }$B$ \textit{be
sets. The difference of }$A$ \textit{and }$B$ \textit{is the set of elements
of }$A$ \textit{that are not in }$B$\textit{, denoted by }$A-B$. \textit{We
write}%
\begin{equation*}
A-B=\left\{ x:x\in A\text{ and }x\notin B\right\} 
\end{equation*}

\section{Formalism}

A question is formally defined as follow;\newline
\textbf{Definition 3.1 (Question):\qquad}\textit{A question is identified
with a topology }$\mathcal{T}_{i}$ \textit{on a given set X. The collection }%
$\mathfrak{T}=\left\{ \mathcal{T}_{i}:i\in I,\text{ }\mathcal{T}_{i}\text{
is a topology on }X\right\} $ \textit{is the set of all possible questions
that one can ask.}

Definition-3.1 requires objetive type questions. For illustration purpose
consider the following objective-type question;\newline
Q1:\qquad What is a particle?\newline
A\qquad a particle is something which has mass.\newline
B\qquad a particle is something which has spin.\newline
C\qquad a particle has both mass and spin.\newline
D\qquad none of the above.\newline
Let mass = $m$, and spin = $s$, then the ground set for the above question
is $X=\left\{ m,s\right\} $. Question Q1 is described by the following
topology 
\begin{equation*}
\mathcal{T}_{1}=\left\{ A,B,C,D\right\} =\left\{ \left\{ m\right\} ,\left\{
s\right\} ,X,\phi \right\} 
\end{equation*}%
Let me ask the other possible questions in this space.\newline
Q2:\qquad What is a massive particle?\newline
A:\qquad having mass.\newline
B:\qquad having mass and spin.\newline
C\qquad none of the above.\newline
Q2 is described by the topology $\mathcal{T}_{2}=\left\{ \left\{ m\right\}
,X,\phi \right\} $. Similarly the question, \textit{'What is a spinning
particle?'}, is described by $\mathcal{T}_{3}=\left\{ \left\{ s\right\}
,X,\phi \right\} $. And the last question, \textit{'What is a
spinning-massive particle?'}, is described by $\mathcal{T}_{4}=\left\{
X,\phi \right\} $.

Every question consists of several issues. An issue is formally defined as
follow;\newline
\textbf{Definition 3.2 (Issue):\qquad }\textit{Any arbitrary point }$x\in X$ 
\textit{be an irreducible assertion and the nbhd }$N$\textit{\ of point }$x$%
\textit{\ is said to be an issue related with point }$x$\textit{\ and the
nbhd system} $\mathcal{N}_{\left( x\right) }$ \textit{is the set of all
issues concerning }$x$. \textit{The open nbhd is said be a relevant issue
for the current problem otherwise it is irrelevant.}

Having defined a question and an issue, one has to resolve the issues of
various points turn by turn by asking various questions. The following
definition resolves an issue;\newline
\textbf{Definition 3.3 (Resolving an issue):\qquad }\textit{Let} $\left( X,%
\mathcal{T}\right) $ \textit{be a topological space}. $\mathcal{T}-\mathcal{N%
}_{(x_{i})}$ \textit{resolves the issue of} $x_{i}\in X$, \textit{where} $%
\mathcal{N}_{(x_{i})}$ \textit{is the neighborhood system of} $x_{i}$.

Definition-3.3 tells us when an issue is solved than that issue no more
exists in the space and hence we are not worried about it any more. One can
see that the elements of $\mathcal{T}-\mathcal{N}_{(x_{i})}$ do not contain $%
x_{i}$ as a point. It is a kind of elimination procedure. The above
definition leads to the following important theorem which may play essential
role in machine learning. But let me first prove a lemma.\newline
\textbf{Lemma 3.1:\qquad }\textit{Let} $X\subseteq Y$\textit{, let} $%
\mathcal{T}_{X}$ \textit{be a topology on} $X$, \textit{then there exists a
topology} $\mathcal{T}_{Y}$ \textit{on} $Y$\textit{, such that} $\mathcal{T}%
_{X}\subseteq \mathcal{T}_{Y}$.\newline
Said differently, every open set of $X$\textit{\ }can be open in $Y$.\newline
\textbf{Proof\qquad }Let $\mathcal{T}_{X}$ and $\mathcal{T}_{Y}$ be discrete
topologies on $X$ and $Y$ respectively. Let $X\subseteq Y$. Since $X\cap Y=X$
($X\subseteq Y$ by hypothesis) also given that $\mathcal{T}_{X}$ and $%
\mathcal{T}_{Y}$ are discrete, therefore $X$ belongs to $\mathcal{T}_{Y}$
and hence every member of $\mathcal{T}_{X}$ is also a member of $\mathcal{T}%
_{Y}$. Hence $\mathcal{T}_{X}\subseteq \mathcal{T}_{Y}$. \newline
\textbf{Theorem 3.1:\qquad }\textit{Let} $\left( X,\mathcal{T}_{X}\right) $ 
\textit{be a topological space,\newline
}i)\qquad \textit{If }$x_{i}\in X$ \textit{belongs to some but not all open
sets of }$X$\textit{\ then }$\mathcal{T}_{X}-\mathcal{N}_{(x_{i})}$\textit{\
is a subspace topology of }$X$\textit{.\newline
}ii)\textit{\qquad If each non-empty open set of }$X$ \textit{contains} $%
x_{i}$ \textit{as a point, then} $\mathcal{T}_{X}-\mathcal{N}_{(x_{i})}$ 
\textit{is the singleton} $\left\{ \phi \right\} $.\newline
iii)\qquad \textit{If there exists no open set that contains }$x_{i}$ 
\textit{as a point}, \textit{then} $\mathcal{T}_{X}-\mathcal{N}_{(x_{i})}$ 
\textit{is empty}.\newline
\textbf{Proof\qquad }i)\qquad Let $U=\left\{ u_{\alpha }:\alpha \in
I\right\} $ and $V=\left\{ v_{\alpha }:\alpha \in I\right\} $ be
subcollections of $\mathcal{T}_{X}$ such that $x_{i}\in u_{\alpha }$ for
every $u_{\alpha }\in U$ and $x_{i}\notin v_{\alpha }$ for every $v_{\alpha
}\in V$. Let $A=\bigcup\limits_{v\in V}v$, to show (a) $\mathcal{T}_{X}-%
\mathcal{N}_{(x_{i})}$ is a topology on $A$, (b) $A$ is subspace of $X$.
First prove part-a\newline
C1:\qquad Let $x_{i}\in u_{\alpha }\in U$ and $x_{i}\notin v_{\alpha }\in V$%
. It implies that $u_{\alpha }$ is a nbhd of $x_{i}$, therefore $u_{\alpha
}\in \mathcal{N}_{(x_{i})}$ and $x_{i}\notin v_{\alpha }$ implies that $%
v_{\alpha }$ is not a nbhd of $x_{i}$, therefore $v_{\alpha }\notin $ $%
\mathcal{N}_{(x_{i})}$. Hence $u_{\alpha }\notin \mathcal{T}_{X}-\mathcal{N}%
_{(x_{i})}$ and $v_{\alpha }\in \mathcal{T}_{X}-\mathcal{N}_{(x_{i})}$.
Hence $A=\bigcup\limits_{v\in V}v\in \mathcal{T}_{X}-\mathcal{N}_{(x_{i})}$.
Since $\mathcal{N}_{\left( x_{i}\right) }$ is non-empty (by theorem-2.1) and 
$\phi $ belongs to $\mathcal{T}_{X}$ (by theorem-2.1(C1)). Thereofore $%
\mathcal{T}_{X}-\mathcal{N}_{(x_{i})}$ also contains $\phi $.\newline
C2:\qquad Given that $x_{i}\notin $ $v_{\alpha }$ , therefore $x_{i}\notin
\bigcup\limits_{\alpha }v_{\alpha }\in \mathcal{T}_{X}-\mathcal{N}_{(x_{i})}$%
.\newline
C3:\qquad Similarly $x_{i}\notin v_{1}\cap v_{2}\in \mathcal{T}_{X}-\mathcal{%
N}_{(x_{i})}$.\newline
Hence $\mathcal{T}_{X}-\mathcal{N}_{(x_{i})}$ is a topology on $A$.\newline
(b)\qquad It is shown in part-a that $\mathcal{T}_{X}-\mathcal{N}_{(x_{i})}$
is a topology on $A$. Further $A\subseteq X$, therefore $\mathcal{T}_{X}-%
\mathcal{N}_{(x_{i})}=\left\{ A\cap U:U\in \mathcal{T}_{X}\right\} $ is a
subspace topology of $X$ (by the definition of subspace topology).\newline
ii)\qquad Let $x_{i}\in u$ for every non-empty set $u\in \mathcal{T}_{X}$,
then by theorem-2.1 each $u$ is contained in $\mathcal{N}_{(x_{i})}$. Since $%
\mathcal{N}_{(x_{i})}$ is non-empty hence $\phi $ is the only member of $%
\mathcal{T}_{X}$ which is not in $\mathcal{N}_{(x_{i})}$. Therefore $%
\mathcal{T}_{X}-\mathcal{N}_{(x_{i})}$ is the singleton $\left\{ \phi
\right\} $.\newline
iii)\qquad Let $x_{i}\notin u$ for every $u\in \mathcal{T}_{X}$ and assume $%
\mathcal{T}_{X}-\mathcal{N}_{(x_{i})}$ is non-empty. When $x_{i}\notin u$,
then there exists no nbhd $N$ \ of $x_{i}$ such that $N\in \mathcal{N}%
_{\left( x_{i}\right) }$. Therefore $\mathcal{N}_{\left( x_{i}\right) }$ is
empty which is a contradiction ( because $\mathcal{N}_{\left( x_{i}\right) }$
is non-empty by theorem-2.1). Hence $\mathcal{T}_{X}-\mathcal{N}_{(x_{i})}$
is empty.\newline
Converse of the theorem is not unique by lemma-3.1 in that $\mathcal{T}_{X}-%
\mathcal{N}_{(x_{i})}\subseteq \mathcal{T}_{X}\subseteq \mathcal{T}_{Y}$.%
\newline

Theorem-3.1 has three parts. It admits three types of questions i.e.,
type-I, type-II and type-III given by cases (i), (ii) and (iii)
respectively. Type-I question generates a sub-question, type-II question has
a definite answer and type-III question is irrelevant. Here if $\mathcal{T}-%
\mathcal{N}_{(x_{i})}$ ends up with sub-question (a sub-question is
identified with a subspace), then one may ask further questions we call them
sub-questions. When $\mathcal{T}-\mathcal{N}_{(x_{i})}$ is equal to $\{\phi
\}$, it means that the question has a definite answer and thus no further
questions are needed to ask and when $\mathcal{T}-\mathcal{N}_{(x_{i})}$ is
empty then the question was irrelevant. To better understand various types
of a question, let me increase the size of the space in the previous
example, by including charge of the particle. Let $X=\left\{ m,s,e\right\} $%
, let $\mathcal{T}_{X}=\left\{ \left\{ m\right\} ,\left\{ m,s\right\}
,\left\{ m,e\right\} ,X,\phi \right\} $ which corresponds to the question; 
\textit{What is a massive particle? }I am interested to resolve the issue of
charge, then%
\begin{align*}
\mathcal{T}_{X}-\mathcal{N}_{(e)}& =\left\{ \left\{ m\right\} ,\left\{
m,s\right\} ,\left\{ m,e\right\} ,X,\phi \right\} -\left\{ \left\{
m,e\right\} ,X\right\} \\
& =\left\{ \left\{ m\right\} ,\left\{ m,s\right\} ,\phi \right\}
\end{align*}%
Since $\mathcal{T}_{X}-\mathcal{N}_{(e)}$ is a topology on $\{m,s\}$,
therefore it corresponds to the sub-question; '\textit{What is a
neutral-massive particle?' One can see that }$\mathcal{T}_{X}-\mathcal{N}%
_{(x_{i})}$ is an elimination procedure, i.e. $\mathcal{T}_{X}-\mathcal{N}%
_{(x_{i})}$ is a collection of subset of $X$ which deos not contain $x_{i}$
as a point. In the above example, when the issue of charge is solved then
the charge no more exits in the subspace. If we remove the nbhd system of $m$%
, then%
\begin{equation*}
\mathcal{T}_{X}-\mathcal{N}_{(m)}=\left\{ \phi \right\} 
\end{equation*}%
Since $\mathcal{T}_{X}-\mathcal{N}_{(m)}$ describes type-II question. Hence
a massive particle has a definite mass and so no further questions are
needed to ask. Now consider type-III question, suppose it is given that the
motion of a particle is described by classical mechanics, then asking about
the quantum mechanical nature of the particle is irrelevant.

The converse of the theorem is not unique means that while starting from a
given solution one can generate several higher order questions.\newline

\subsection{Most and least efficient questions}

The operation $\mathcal{T}-\mathcal{N}_{(x_{i})}$ acts most efficiently and
terminates the processes in just one single step, when theorem-3.1(ii) is
true. In the collection $\mathfrak{T}=\left\{ \mathcal{T}_{j}:j\in I\right\} 
$, there exist $\mathcal{T}_{k}$ such that $\mathcal{T}_{k}-\mathcal{N}%
_{(x_{i})}=\left\{ \phi \right\} $. If $X$ is discrete, then $\mathcal{T}-%
\mathcal{N}_{(x_{i})}$ operates least efficiently. When $X$ has more then
one points then discrete topologies satisfy the first case of theorem-3.1.
When $X$ is a singleton then it satisfies both first and second cases of
theorem-3.1. One can observe that when $X$ is discrete, then the subspace
topology $\mathcal{T}-\mathcal{N}_{(x_{i})}$ is a topology on $X-\left\{
x_{i}\right\} $. Thus $\mathcal{T}-\mathcal{N}_{(x_{i})}$ eliminates just
one irreducible assertion from the space. There is a topological property,
called the hereditary property, whenever a topological space $X$ has
property $P$, then every subspace of $X$ also has property $P$. It implies
that $\mathcal{T}-\mathcal{N}_{(x_{i})}$ is a discrete topology on $%
X-\{x_{i}\}$. Let $\mathcal{T}$ be a discrete topology on $X$, define $%
\mathcal{T}_{2}=\mathcal{T}-\mathcal{N}_{(x_{i})}$, then $\mathcal{T}_{2}-%
\mathcal{N}_{(x_{i})}$ is a discrete topology on $X-\{x_{i}\}$, such a
process is very slow. Every time just one irreducible assertion is
eliminated from the space.

\subsection{Negation question}

A question can also be asked negatively. Consider two systems $A$ and $B$
which are thermally connected. Let the energy of $A$ is greater than the
energy of $B$\textit{\ }(i.e., $E_{A}>E_{B}$), then the energy is flowing
from $A$ to $B$. In this case system $A$ asks from system $B$, "Do you want
to gain energy?" In the response system \textit{B} asks from system \textit{A%
}, "Don't you want to gain energy?"\textit{\ }If we define the former
question is to be a topology $\mathcal{T}_{X}^{U}$, then the later question
is defined as follow;\newline
\textbf{Definition 3.4 (negation question):\qquad }\textit{Let }$\mathcal{T}%
_{X}^{U}$ \textit{be a topology on} $X$. \textit{The corresponding negation
question is defined as }$\mathcal{T}_{X}^{F}=\left\{ X-U:U\in \mathcal{T}%
_{X}^{U}\right\} $.

One can show that $\mathcal{T}_{X}^{F}$ is also a topology on $X$ and that $%
\mathcal{T}_{X}^{U}\neq \mathcal{T}_{X}^{F}$ in general. The intersection $%
\mathcal{T}_{X}^{U}\cap \mathcal{T}_{X}^{F}$ is the collection of clopen
sets. When every member of $\mathcal{T}_{X}^{U}$ is clopen, then$\mathcal{T}%
_{X}^{U}=\mathcal{T}_{X}^{F}$. In this case question and its negation
question are identical. It is believed that question and negation question
are always identical \cite{Knuth}. It is argued that both questions are
identical in the sense that they ask the same thing. It is demonstrated by
asking, '\textit{Is it raining?'} and \textit{'Is it not raining?'} Both
questions are answered by the statements \textit{'It is raining!'} and 
\textit{'It is not raining!'}, and thus they are equivalent. Here an
important question is missing, that is, \textit{'Who asks the question?'}.
Let me replace the raining question by the following simple question and
demonstrate that they are not equivalent. Consider Bob and Alice. Let Bob
has extra money and Alice needs some money. Bob wants to get rid of extra
money and asks Alice, "Do you need money?" \ Alice needs the money she would
respond, "Don't you need the money?" The two persons ask questions according
to their desires and thus generate different questions. Here money is
flowing from one person to the other. We conclude from this that for a
machine (Machine is defined below) to be intelligent it must have desires. A
clear example of an intelligent machine who has got desires is human. Let me
put it in this way humans are intelligent because they have desires. A truly
artificially intelligent machine can think independently when it is able to
create desires which forces it to act accordingly. We also claim that
desires are quantifiable.

The equivalence of question and negation question make sense in probability
theory. A probability space is a triple $\left( \Omega ,\mathcal{F},P\right) 
$, where $\mathcal{F}$ is a collection of subsets of $\Omega $ relevant to a
particular experiment and $P$ is a probability measure. The elements of $%
\mathcal{F}$ are called events. A non-empty collection of subsets $\mathcal{F%
}$ is called a $\sigma $-field on $\Omega $ if \cite{Grigoriu}\newline
1.\qquad $\phi \in \mathcal{F}$\newline
2.\qquad $A\in \mathcal{F}$, then $A^{c}\in \mathcal{F}$, where $A^{c}$ is
compliment of $A$.\newline
3.\qquad $A_{i}\in \mathcal{F}$, $i\in I$, $I=$ a countable set, then $%
\bigcup\limits_{i\in I}A_{i}\in \mathcal{F}$.

The last condition implies $\left( \bigcup\limits_{i\in I}A_{i}\right)
^{c}=\tbigcap\limits_{i\in I}A_{i}^{c}$ (by De Morgan's formula). Therefore $%
\tbigcap\limits_{i\in I}A_{i}^{c}\in \mathcal{F}$. In our case it
corresponds to $\mathcal{T}_{X}^{U}=\mathcal{T}_{X}^{F}$.

\subsection{Machine}

The negation question is very important for construction of a machine. A
machine is defined as an agent who asks question $\mathcal{T}_{X}^{U}$. For
every machine there exists an anti-machine who asks the negation question $%
\mathcal{T}_{X}^{F}$. A machine and an anti-machine make a universe. Here
'machine' is used in a very broad sense. It can be an intelligent machine
which ask a question that has definite answer (type-II questions) or it can
be a system of particles and anti-machine is the rest of the universe or the
environment. We suspect that noise is also a kind of anti-machine. The two
machines communicates through the clopen sets. Said differently the clopen
sets are the information the two machines share to operate and the
non-clopen set are the information they do not share. The two machine are in
perfect agreement when they share all informations. In this case $\mathcal{T}%
_{X}^{U}=\mathcal{T}_{X}^{F}$.

The collection\textit{\ }$\mathfrak{T}=\left\{ \mathcal{T}_{i}:i\in
I\right\} $ is the set of all possible questions that a machine can ask.
Therefore a machine consists of several components (subsystems). Each
component is called an atomic machine. Each atomic machine asks a question $%
\mathcal{T}_{i}^{U}\in \mathfrak{T}$, and for every atomic machine there
exists an anti-atomic machine who asks the negation question $\mathcal{T}%
_{i}^{F}\in \mathfrak{T}$. Some atomic machine are their own anti-atomic
machine for which $\mathcal{T}_{j}^{U}=\mathcal{T}_{j}^{F}\in \mathfrak{T}$.

\section{Conclusion}

It is possible to construct a machine that asks intelligent questions and
work most efficiently. We suggest a machine work most efficiently if it asks
type-II question. This paper is a qualitative development of question
theory. Our formalism is general, it can be applied to any problem where
question exists.\newline

\textbf{Acknowledgements}

I am indebted to A. Caticha, A. Inomata and S. Ali for many insightful
remarks and valuable suggestions. I also benefited from K. Knuth's
thoughtful lectures on Question Algebra.

\end{document}